\documentclass[12pt,leqno]{article}
\usepackage{amsmath,amssymb}
\begin{document}
\newtheorem{theorem}{Theorem}[section]
\newtheorem{lemma}{Lemma}[section]
\newtheorem{corol}{Corollary}[section]
\newtheorem{prop}{Proposition}[section]
\newtheorem{defin}{Definition}[section]
\newtheorem{rem}{Remark}[section]
\newtheorem{example}{Example}[section]
\title{Dirac submanifolds of Jacobi manifolds}
\author{{\small by}\vspace{2mm}\\Izu Vaisman}
\date{}
\maketitle
{\def\thefootnote{*}\footnotetext[1]%
{{\it 2000 Mathematics Subject Classification: 53D17} .
\newline\indent{\it Key words and phrases}: Soldered tensor field.
Poisson manifold. Jacobi manifold.
Dirac submanifold.}}
\begin{center} \begin{minipage}{12cm}
A{\footnotesize BSTRACT. The notion of a Dirac submanifold of a
Poisson manifold studied by Xu is interpreted in terms of a
general notion of tensor fields soldered to a normalized
submanifold. This interpretation is used to define the notion of a
Dirac submanifold of a Jacobi manifold. Several properties and
examples are discussed.}
\end{minipage}
\end{center} \vspace{5mm}
\section{Normalized submanifolds}
In this section we make some general considerations on
submanifolds of a differentiable\footnote{In this paper everything
is of class $C^\infty$ and all submanifolds are embedded.}
manifold. These considerations were inspired by the theory of
Dirac submanifolds of a Poisson manifold developed in \cite{Xu}.
The notion of a normalized submanifold was used in affine and
projective differential geometry half a century ago.
\begin{defin} \label{normalized} {\rm Let $N^n$ be a submanifold
of $M^m$ (indices denote dimensions), and $\iota:N\subseteq M$ the
corresponding embedding. A {\it normalization} of $N$ by a {\it
normal bundle} $\nu N$ is a splitting
\begin{equation} \label{split} TM|_N=TN\oplus\nu N
\end{equation} ($T$ denotes tangent bundles). A submanifold
endowed with a normalization is called a {\it normalized
submanifold}}.
\end{defin}

Then, if $X\in\Gamma TM$ ($\Gamma$ denotes spaces of global cross
sections) is a vector field on $M$ such that $X|_N\in\Gamma TN$,
respectively $X|_N\in\Gamma \nu N$, $X$ is said to be {\it
tangent}, respectively {\it normal}, to $N$.

Let $(N^n,\nu N)$ be a normalized submanifold of $M^m$. Let
$\sigma:W\rightarrow N$ be a tubular neighborhood of $N$ such that
$\forall x\in N$, $T_x(W_x)=\nu_x N$ ($W_x$ is the fiber of $W$
and $\nu_xN$ is the fiber of $\nu N$ at $x$). Then $W$ is said to
be a {\it compatible tubular neighborhood} of $N$; obviously, such
neighborhoods exist. Furthermore, each point $x\in N$ has a
$\sigma$-trivializing neighborhood $U$ endowed with coordinates
$(x^a)$ $(a,b,c,...= 1,...,m-n)$ on the fibers of $\sigma$ and
such that $x^a|_{N\cap U}=0$, and coordinates $(y^u)$
$(u,v,w,...=m-n+1,...,m)$ on $N\cap U$. We will say that
$(x^a,y^u)$ are {\it adapted local coordinates}.

With respect to adapted coordinates, $N$ has the local equations
$x^a=0$, and
\begin{equation} \label{adapted} TN|_{U\cap N}=span\left\{\left.\frac{\partial}
{\partial y^u}\right|_{x^a=0}\right\},\;\; \nu N|_{U\cap N}
=span\left\{\left.\frac{\partial} {\partial
x^a}\right|_{x^a=0}\right\}.\end{equation} Accordingly, the
transition functions between two systems of adapted local
coordinates must be of the form \begin{equation}
\label{transition} \tilde x^a=\tilde x^a(x^b,y^v),\;\tilde
y^u=\tilde y^u(y^v),
\end{equation}  and satisfy the conditions \begin{equation} \label{transcond}
\left.\frac{\partial\tilde x^a}{\partial y^v}\right|_{x^b=0}=0,
\;\frac{\partial\tilde y^u}{\partial x^b}\equiv0.
\end{equation}

The splitting (\ref{split}) induces a similar relation for the
dual bundles: \begin{equation} \label{dualsplit}
T^*M|_N=T^*N\oplus\nu^*N,\end{equation} and, locally, with respect
to adapted coordinates one has
\begin{equation}\label{dualspan}
\begin{array}{c} T^*N=ann(\nu N)= span\{dy^u|_{x^a=0}\},\vspace{2mm}\\
\nu^*N=ann(TN)=span\{dx^a|_{x^a=0}\} \end{array}
\end{equation} ($ann$ denotes annihilator spaces).

Two normal bundles $\nu N$ and $\tilde\nu N$ of the same
submanifold $N$ are connected as follows. Let $p_\nu,p_T$ be the
projections defined by the splitting (\ref{split}), and
$p_{\tilde\nu},\tilde p_T$ the similar projections of the second
normalization. The mapping $(v\in\nu N)\mapsto p_{\tilde\nu}v$ is
an isomorphism $\varphi:\nu N\rightarrow\tilde\nu N$, and with
respect to adapted coordinates we have
\begin{equation}\label{nutilde} \tilde\nu
N=span\{\varphi\left(\left.\frac{\partial}{\partial
x^a}\right|_{x^c=0}\right)= X_a|_{x^c=0}\},\end{equation}
where\footnote{In this paper, we use the Einstein summation
convention.}
\begin{equation} \label{Xmare} X_a=\frac{\partial}{\partial
x^a}-\theta^u_a\frac{\partial}{\partial y^u}\end{equation} are
vector fields on $U$, and $(\theta^u_a|_{x^c=0})$ is the local
matrix of the homomorphism $\psi:\nu N\rightarrow TN$ defined by
$\psi(v)=\tilde p_T(v)$ $(v\in\nu N)$. Notice that $\varphi+\psi$
is the inclusion of $\nu N$ in $TM|_N$, and $\tilde\nu N$ is
uniquely determined by any of the mappings $\varphi,\psi$.

Some geometric objects of $M$ may have a strong relationship with
the normalized submanifold $(N,\nu N)$.
\begin{defin}\label{formsold} {\rm A differential $k$-form
$\kappa\in\Omega^k(M)$ ($\Omega$ denotes spaces of differential
forms) is {\it soldered} to $N$ if $\iota^*[L_X\kappa]=0$ for any
vector field $X\in\Gamma TM$ normal to $N$.}
\end{defin}

Since $\forall f\in C^\infty(M)$ and any vector fields
$X,Y_1,\ldots,Y_k$, one has
\begin{equation}\label{LiefX}
(L_{fX}\kappa)(Y_1,...,Y_k)=f(L_{X}\kappa)(Y_1,...,Y_k)
\end{equation}
$$-\sum_{j=1}^k(-1)^j(Y_jf)[i(X)\kappa](Y_1,...,Y_{j-1},Y_{j+1},
...,Y_k),$$ it follows that $\kappa$ is soldered iff for any
vector field $X\in\Gamma TM$ normal to $N$ one has
\begin{equation} \label{formekappa}
\iota^*[i(X)\kappa]=0,\;\iota^*[L_X\kappa]=0, \end{equation} or,
equivalently, \begin{equation} \label{equivalforme}
\iota^*[i(X)\kappa]=0,\;\iota^*[i(X)d\kappa]=0.
\end{equation}

With respect to adapted local coordinates, $\kappa$ has the
expression
\begin{equation}\label{kappalocal} \kappa=\sum_{s+t=k}\frac{1}{s!t!}
\kappa_{a_1...a_su_1...u_t}dx^{a_1}\wedge\cdots\wedge dx^{a_s}
\wedge dy^{u_1}\wedge\cdots\wedge dy^{u_t},\end{equation} and
$\kappa$ is soldered to $N$ iff \begin{equation}
\label{soldkappalocal} \kappa_{au_1...u_{k-1}}|_{x^b=0}=0,\;
\left.\frac{\partial\kappa_{u_1...u_k}}{\partial
x^a}\right|_{x^b=0}=0.\end{equation}

In particular, the space of soldered functions is
\begin{equation} \label{soldfunct} C^\infty(M,N,\nu N)=
\{ f\in C^\infty(M)\;/\;(\partial f/\partial x^a)_{x^a=0}=0\}.
\end{equation}

We will denote by $\Omega^k(M,N,\nu N)$ the space of soldered
$k$-forms (for $k=0$ we have the space (\ref{soldfunct})).
Obviously, the soldering conditions (\ref{equivalforme}) are
compatible with the exterior product and the exterior
differential. Therefore, we get a cohomology algebra
$\oplus_kH^k_{sdeR}(M,N,\nu N)$, which will be called the {\it
soldered de Rham cohomology algebra}, defined by the cochain
spaces $\Omega^k(M,N,\nu N)$ and the operator $d$. The inclusion
in the usual de Rham complex induces homomorphisms
\begin{equation} \label{homdeR} \iota^k:H^k_{sdeR}(M,N,\nu
N)\longrightarrow H^k_{deR}(M).\end{equation}

In principle, the spaces $H^k_{sdeR}(W,N,\nu N)$, where $W$ is a
compatible tubular neighborhood, should provide interesting
information about the normalized submanifold $(N,\nu N)$.
\begin{defin} \label{vectsold}
{\rm A $k$-vector field $Q\in\mathcal{V}^k(M)$ ($ \mathcal{V}$
denotes spaces of multivector fields) is {\it soldered} to the
normalized submanifold $(N,\nu N)$ if: i) for any $(k-1)$-form
$\lambda\in\Gamma\wedge^{k-1}[ann(\nu N)]$ the vector field
defined along $N$ by $i(\lambda)Q|_N$ is tangent to $N$; ii) for
any vector field $X$ on $M$ normal to $N$, $(L_XQ)|_{ann(\nu
N)}=0$.}\end{defin}

We will denote by $ \mathcal{V}^k(M,N,\nu N)$ the space of
soldered $k$-vector fields.  Using adapted local coordinates we
see that $ \mathcal{V}^0(M,N,\nu N)$ is, again, (\ref{soldfunct}),
and for $k=1$ one has
\begin{equation} \label{Diractens1}
\mathcal{V}^1(M,N,\nu N )=\left\{Y\in\mathcal{V}^1(M)\;/\;Y|_N
=\eta^u\left.\frac{\partial}{\partial
y^u}\right|_N,\,\left.\frac{\partial\eta^u}{\partial
x^a}\right|_{x^c=0}=0\right\}. \end{equation}

Generally, we have an expression of the form
\begin{equation} \label{Qlocal} Q=\sum_{s+t=k}\frac{1}{s!t!}
Q^{a_1...a_su_1...u_t}\frac{\partial}{\partial x^{a_1}}
\wedge\cdots\wedge\frac{\partial}{\partial x^{a_s}} \wedge
\frac{\partial}{\partial y^{u_1}}
\wedge\cdots\wedge\frac{\partial}{\partial y^{u_t}},\end{equation}
and $Q\in\mathcal{V}^k(M,N,\nu N)$ iff
\begin{equation} \label{conditionsQ}
Q^{au_1...u_{k-1}}|_{x^c=0}=0,\;\left.\frac{\partial
Q^{u_1...u_k}}{\partial x^a}\right|_{x^c=0}=0.\end{equation}

The spaces of soldered forms and multivector fields are components
of important algebraic structures namely,
\begin{prop} \label{propos00} 1. The space $\mathcal{V}^1(M,N,\nu N )$
is a Herz-Reinhart Lie algebra over $( \mathbb{R},
C^\infty(M,N,\nu N))$. 2. The complex of $N$-soldered differential
forms is a complex over the Lie algebra $\mathcal{V}^1(M,N,\nu N
)$. 3. The triple $(\oplus_k\mathcal{V}^k(M,N,\nu
N),\wedge,[\,,\,])$, where $[\;,\;]$ is the Schouten-Nijenhuis
bracket is a Gerstenhaber algebra. \end{prop} \noindent{\bf
Proof.} For the definition of the algebraic structures above see,
for instance, \cite{Dorf} and \cite{Hu}. The use of adapted local
coordinates shows that if $f\in C^\infty(M,N,\nu N)$ and $Y,Z\in
\mathcal{V}^1(M,N,\nu N)$ then $fY$ and $[Y,Z]$ belong to
$\mathcal{V}^1(M,N,\nu N)$, which proves 1. Furthermore, for the
same $Y$ and any $\kappa\in\Omega^k(M,N,\nu N)$, $i(Y)\kappa\in
\Omega^{k-1}(M,N,\nu N)$, and we get 2. Finally, the exterior
product of soldered multivector fields obviously is soldered, and
in order to get 3, it remains to prove that the Schouten-Nijenhuis
bracket\footnote{We take this bracket with the sign convention of
the axioms of graded Lie algebras, e.g., \cite{V1}, Proposition
4.21.} of $P\in\mathcal{V}^p(M,N,\nu N)$ and
$Q\in\mathcal{V}^q(M,N,\nu N)$ belongs to $
\mathcal{V}^{p+q-1}(M,N,\nu N)$. To see that soldering condition
i) is satisfied, we look at the known formula (e.g., \cite{V1})
\begin{equation} \label{Lichne} i([P,Q])\omega=
(-1)^{(p+1)(q+1)}i(P)d[i(Q)\omega]-i(Q)d[i(P)\omega]\end{equation}
$$+(-1)^pi(P\wedge Q)d\omega,$$ where $\omega$ is an arbitrary
$(p+q-1)$-form on $M$, and use this formula for
$$\omega=dx^a\wedge dy^{u_1}\wedge\cdots \wedge dy^{u_{p+q-2}}.$$
Then, soldering condition ii) follows by using (\ref{Lichne}) to
evaluate the terms of the equality
$$L_X[P,Q]=[L_XP,Q] + [P,L_XQ],$$ where $X$ is normal
to $N$, on $dy^{u_1}\wedge\cdots\wedge dy^{u_{p+q-1}}$. Q.e.d.

The following proposition extends a result proven for Poisson
bivector fields in \cite{Xu} \begin{prop} \label{proposinvol} If
the involutive diffeomorphism $\varphi:M\rightarrow M$
$(\varphi^2=Id)$ preserves a $k$-form $\kappa$, respectively a
$k$-vector field $Q$, then $\kappa$, respectively $Q$, is soldered
to the fixed point locus $N$ of $\varphi$. \end{prop}
\noindent{\bf Proof.} It is well known that $N$ is a submanifold
of $M$, the tangent bundle $TN$ consists of the $(+1)$-eigenspaces
of $\varphi_*$ along $N$, and $N$ has a normalization with the
normal bundle $\nu N$ defined by the $(-1)$-eigenspaces of
$\varphi_*$ along $N$ \cite{Xu}. This also implies that $ann(\nu
N)=T^*N$ consists of the $(+1)$-eigenspaces, and $ann(TN)=\nu^*N$
consists of the $(-1)$-eigenspaces, of $\varphi^*$ along $N$ . If
$\varphi^*\kappa=\kappa$, if $X$ is a normal vector field of $N$
on $M$ and $Y_1,...,Y_{k-1}$ are tangent to $N$, then
$$\kappa(X,Y_1,...,Y_{k-1})|_N=(\varphi^*\kappa)(X,Y_1,...,Y_{k-1})|_N
=-\kappa((X,Y_1,...,Y_{k-1})|_N$$ hence $\iota^*[i(X)\kappa]=0$.
The same holds for $d\kappa$, therefore, $\kappa$ is soldered to
$N$. The proof for a $k$-vector field $Q$ is similar, and uses the
fact that
$\varphi_*(L_XQ)=L_{\varphi_*X}(\varphi_*Q)\circ\varphi$, for any
diffeomorphism $\varphi:M\rightarrow N$. Q.e.d.

The $N$-soldered differential forms and multivector fields have a
nice interpretation by means of the geometry of the tangent bundle
$TM$, and by looking at the normal bundle $\nu N$ as a submanifold
of the former. (See \cite{Xu} for the case of a Poisson bivector
field.)

The tensor fields of the manifold $M$ may be lifted to $TM$ by
various processes and, in particular, there exists a {\it complete
lift}, which comes from the lift of the flow of a vector field
\cite{YI}. In the case of differential forms and multivector
fields, the complete lift has the following coordinate expression
\begin{equation} \label{completeform}
\kappa^C=\frac{1}{(k-1)!}\kappa_{i_1...i_k}d\hspace{-1mm}\stackrel{.}{z}^{i_1}
\wedge dz^{i_2}\wedge\cdots\wedge dz^{i_k}\end{equation}
$$+\frac{1}{k!}\stackrel{.}{z}^k\frac{\partial\kappa_{i_1...i_k}}
{\partial z^k}dz^{i_1}\wedge\cdots\wedge dz^{i_k},$$
\begin{equation} \label{completemult} Q^C=\frac{1}{(k-1)!}
Q^{i_1...i_k}\frac{\partial}{\partial z^{i_1}}\wedge
\frac{\partial} {\partial\hspace{-1mm}\stackrel{.}{ z}^{i_2}}
\wedge\cdots\wedge\frac{\partial}
{\partial\hspace{-1mm}\stackrel{.}{ z}^{i_k}} \end{equation}
$$+\frac{1}{k!}\stackrel{.}{z}^k\frac{\partial Q^{i_1...i_k}}
{\partial z^k} \frac{\partial} {\partial\hspace{-1mm}\stackrel{.}{
z}^{i_1}} \wedge\cdots\wedge\frac{\partial}
{\partial\hspace{-1mm}\stackrel{.}{ z}^{i_k}}.$$ In formulas
(\ref{completeform}), (\ref{completemult}) $\kappa$ and $Q$ are
given by the expressions (\ref{kappalocal}) and (\ref{Qlocal}),
respectively, while $(z^i)=(x^a,y^u)$ $(i=1,...,m)$, and
$(\stackrel{.}{z}^i)$ are the corresponding natural vector
coordinates. \begin{prop} \label{propos01} i) The differential
form $\kappa$ is soldered to the normalized submanifold $(N,\nu N)
$ of $M$ iff, $\forall Z\in\Gamma T(\nu N)$, the form $i(Z)
\kappa^C$ belongs to the ideal generated by $\Gamma[ann( T(\nu
N))]$. ii) The $k$-vector field $Q$ is soldered to $(N,\nu N)$
iff, $\forall\alpha\in\Gamma [ann(T(\nu N))]$, $i(\alpha)Q^C$
belongs to the ideal generated by $\Gamma T(\nu N)$.
\end{prop} \noindent{\bf Proof.} On $M$, we use $N$-adapted local
coordinates and on $TM$ the corresponding natural vector
coordinates as described above. Then, the submanifold $\nu N$ has
the local equations $x^a=0,\stackrel{.}{y}^u=0$, and the results
are immediate consequences of formulas (\ref{completeform}) and
(\ref{completemult}). Q.e.d.

We remark that it is possible to define soldered symmetric tensor
fields, similarly. The following proposition provides a nice
example:
\begin{prop} \label{propos02} Let $M$ be a Riemannian manifold
with the metric tensor $g$. Then $g$ is soldered to a submanifold
$N$ normalized by the normal bundle $\nu N\perp TN$ iff $N$ is a
totally geodesic submanifold. \end{prop} \noindent{\bf Proof.} If
$\nabla$ is the Levi-Civita connection of $g$, i.e., $\nabla g=0$
and $\nabla$ has no torsion, one has
$$(L_Xg)(Y,Z)=g(\nabla_YX,Z)+g(Y,\nabla_ZX).$$ Then, if $X$ is normal
and $Y,Z$ are tangent to $N$, the restriction of the previous
formula to $N$ yields
\begin{equation} \label{secondf}
(L_Xg)(Y,Z)=2(g(b(Y,Z),X),\end{equation} where $b$ is the second
fundamental form of the submanifold $N$. Thus, the soldering
condition of $g$ is equivalent to $b=0$. Q.e.d.

The notion of soldering has other interesting extensions too.
First, for a multivector field $Q\in \mathcal{V}^k(M)$, condition
ii) of Definition \ref{vectsold} is itself a geometric condition,
since it is easy to check that if it holds for $X$ normal to $N$
it also holds for $fX$, $\forall f\in C^\infty(M)$. If $Q$
satisfies only this condition, we will call it {\it
quasi-soldered} to $(N,\nu N)$. In fact, this notion extends to
any contravariant tensor field. The notion of a soldered
differential form also extends to any covariant tensor field but,
it implies algebraic conditions like the first condition
(\ref{formekappa}) too.

Then, we may look at objects that only satisfy the algebraic
condition of soldering (e.g., the first condition
(\ref{formekappa}), condition i) of Definition \ref{vectsold},
etc.), and call them {\it algebraically compatible with the
normalized submanifold}. If the algebraic condition holds, the Lie
derivative in the normal directions yields an important object for
the submanifold.  For instance, a Riemannian metric is
algebraically compatible with any submanifold $N$, if the normal
bundle is the $g$-orthogonal bundle of $TN$, and formula
(\ref{secondf}) shows that $(L_Xg)|_{TN}$ is equivalent with the
second fundamental form of the submanifold.

Finally, we may add to the algebraic condition a condition that is
weaker than soldering. For instance, a tensor field that is
algebraically compatible with a normalized submanifold will be
called {\it weakly soldered} if, along the submanifold, the Lie
derivative in the normal directions are proportional to the
pullback of the tensor field to the submanifold. For instance, the
Riemannian metric $g$ is weakly soldered to a submanifold $N$ with
the $g$-normal bundle if, for some $1$-form $\alpha$, one has
$(L_Xg)|_{TN}=\alpha(X)g|_{TN}$, for any normal vector field $X$,
and this happens iff $N$ is a totally umbilical submanifold.
\section{Dirac submanifolds of Poisson manifolds}
In this section we recall the definition and characteristic
properties of the Dirac submanifolds of a Poisson manifold studied
by Xu \cite{Xu}, and give a few additional facts. A Dirac
submanifold inherits an induced Poisson structure, and the
cotangent Lie algebroid of the latter may be seen as a Lie
subalgebroid of the cotangent Lie algebroid of the original
manifold. We refer to \cite{V1} for generalities of Poisson
geometry.
\begin{defin} \label{Diracsub}
{\rm A submanifold $N^n$ of a Poisson manifold $M^m$, with the
Poisson bivector field $\Pi$, is a {\it Dirac submanifold} if $N$
has a normalization (\ref{split})
with the following properties:\\
i) $\sharp_{\Pi}(ann(\nu N))\subseteq TN$ ($\sharp_{\Pi}$ is the
morphism $T^*M\rightarrow TM$ defined by $\Pi$); if this condition
holds, we will say that $\nu N$ is {\it  algebraically $\Pi$-compatible};\\
ii) $\forall x\in N$ there exists an open neighborhood $U$ of $x$
in $M$ such that $\forall f,g\in C^\infty(U)$ which satisfy the
conditions $df|_{\nu N}=0,dg|_{\nu N}=0$ one has $d\{f,g\}|_{\nu
N}=0$. (\{f,g\} is the Poisson bracket defined by $\Pi$.)}
\end{defin}
\begin{prop}\label{propos11} {\rm\cite{Xu}} With the notation of Definition
\ref{Diracsub}, the submanifold $N$ is a Dirac submanifold iff
there exists a normalization {\rm (\ref{split})} such that the
Poisson bivector field $\Pi$ is soldered to $(N,\nu N)$.
\end{prop} \noindent{\bf Proof.} $\Pi$ is soldered to $(N,\nu
N)$ iff with the notation and the adapted coordinates as defined
in Section 1, one has
\begin{equation} \label{Pilocal}
\Pi=\frac{1}{2}\Pi^{ab}\frac{\partial}{\partial
x^a}\wedge\frac{\partial}{\partial x^b} +
\Pi^{au}\frac{\partial}{\partial
x^a}\wedge\frac{\partial}{\partial y^u} +
\frac{1}{2}\Pi^{uv}\frac{\partial}{\partial
y^u}\wedge\frac{\partial}{\partial y^v}, \end{equation} where
\begin{equation} \label{Piauzero} \Pi^{au}|_{x^c=0}=0, \end{equation}
\begin{equation} \label{derivPi} \left.\frac{\partial\Pi^{uv}}{\partial
x^a}\right|_{x^c=0}=0. \end{equation}

On the other hand, if $\Pi$ is given by (\ref{Pilocal}), condition
i) of Definition \ref{Diracsub} is equivalent to
$\sharp_{\Pi}(dy^u)\in span\{\partial/\partial y^v\}$ along $N$,
which means (\ref{Piauzero}), and condition ii) is equivalent with
the fact that $d\{y^u,y^v\}|_{x^c=0}=0$, which is (\ref{derivPi}).
Q.e.d.

Furthermore, from the definition of soldered multivector fields we
also get \begin{prop}\label{propos14} {\rm \cite{Xu}} The
submanifold $N$ of $(M,\Pi)$ is a Dirac submanifold iff $N$ has a
$\Pi$-compatible normal bundle $\nu N$ such that, for any, normal
to $N$, vector field $X$ of $M$, one has
\begin{equation}\label{Liecond} (L_X\Pi)|_{ann(\nu N)}=0.
\end{equation}\end{prop}
\begin{corol}\label{propos16} Let
$\{X_\alpha\}$ be a family of Poisson infinitesimal automorphisms
of $(M,\Pi)$. Then, any submanifold $N$ such that $
span\{X_\alpha|_N\}$ is a $\Pi$-compatible normal bundle of $N$ is
a Dirac submanifold.\end{corol} \noindent{\bf Proof.} The
hypotheses of Corollary \ref{propos16} imply the characteristic
conditions stated by Proposition \ref{propos14}. Q.e.d.

The $\Pi$-compatibility hypothesis of Corollary \ref{propos16}
also has the following meaning. A family $\{X_\alpha\}$ of Poisson
infinitesimal automorphisms has an {\it associated generalized
distribution} $D(X_\alpha)$ spanned by the $\Pi$-hamiltonian
vector fields of the functions $f\in C^\infty(M)$ that are
constant along the orbits of the vector fields $X_\alpha$, and
this distribution is involutive. $ span\{X_\alpha|_N\}$ is a
$\Pi$-compatible normal bundle of $N$ iff $D(X_\alpha)\subseteq
TN$. If the family $\{X_\alpha\}$ reduces to one hamiltonian
vector field $X^{\Pi}_h$, we have no submanifolds $N$ as in
Corollary \ref{propos16} since $N$ should be both tangent and
normal to $X^{\Pi}_h$. But, we may have the required type of
submanifolds if the family consists of a single infinitesimal
automorphism $X$ that is not a hamiltonian vector field. For
instance, if $D(X)$ is a regular distribution, it must be a
foliation and the leaves of this foliation are Dirac hypersurfaces
of $M$.

Before going on with the discussion of Dirac submanifold, let us
also consider some of the situations mentioned at the end of
Section 1. \begin{defin} \label{algebr} {\rm A normalized
submanifold $(N,\nu N)$ of the Poisson manifold $(M,\Pi)$ will be
called an {\it algebraically Poisson-compatible (a.P.c.)
submanifold}, respectively a {\it quasi-Dirac submanifold}, if the
Poisson bivector field $\Pi$ is algebraically compatible,
respectively quasi-soldered, to the submanifold.}
\end{defin}

Thus, the a.P.c. property is characterized by i) of Definition
\ref{Diracsub}, respectively, by the local condition
(\ref{Piauzero}), and the quasi-Dirac property is characterized by
the local condition (\ref{derivPi}).

The a.P.c. and Dirac properties of a submanifold may hold for more
than one normal bundle \cite{Xu}. A second normal bundle
$\tilde\nu N$ may be defined by (\ref{nutilde}), and the
corresponding local expression of the Poisson bivector field is
obtained by switching to the bases $(X_a,\partial/\partial y^u)$
in (\ref{Pilocal}). Accordingly, for $\tilde\nu N$, the a.P.c.
condition is equivalent to
\begin{equation} \label{cond1theta}
(\Pi^{ab}\theta_b^u)|_N=0\;\Leftrightarrow
\;\psi\circ(\sharp_{\Pi}|_{ann(TN)})=0,\end{equation} and
condition ii) of Definition \ref{Diracsub} is equivalent to
\begin{equation}\label{cond2theta}
X_c[\Pi^{ab} \theta^u_a\theta^v_b-\Pi^{au}\theta^v_a+
\Pi^{av}\theta^u_a+\Pi^{uv}]|_{x^c=0}=0.\end{equation} In view of
(\ref{Piauzero}), (\ref{derivPi}), and (\ref{cond1theta}),
(\ref{cond2theta}) becomes
\begin{equation} \label{cond3theta}
\left(\frac{\partial\Pi^{ab}}{\partial x^c} \theta^u_a\theta^v_b +
\frac{\partial\Pi^{av}}{\partial x^c} \theta^u_b -
\frac{\partial\Pi^{ua}}{\partial x^c} \theta^v_a -
\frac{\partial\Pi^{uv}}{\partial y^w} \theta^w_c\right)_{x^c=0}=0.
\end{equation}

Condition (\ref{cond1theta}) shows that, if
$\sharp_{\Pi}|_{ann(TN)}$ is a surjection, therefore, an
isomorphism onto $\nu N$ (equivalently, $det(\Pi^{ab})\neq0$),
then $\nu N= \sharp_\Pi(ann(\nu N))$ provides the only
normalization which makes $N$ an a.P.c. submanifold of $M$. The
submanifolds $N$ such that $\sharp_\Pi(ann(\nu N))$ is a
complement of $TN$ in $TM|_N$ are called {\it cosymplectic
submanifolds}, and it is known that they are Dirac submanifolds
\cite{Xu}. Indeed, the a.P.c. property follows from the skew
symmetry of $\Pi$, and (\ref{derivPi}) follows from the following
component of the Poisson condition $[\Pi,\Pi]=0$ in the local
adapted coordinates of (\ref{Pilocal}) \cite{V1}:
\begin{equation} \label{auv} [\Pi,\Pi]^{auv}|_{x^c=0}= 2\left(
\Pi^{ab}\frac{\partial\Pi^{uv}}{\partial x^b}\right)_{x^c=0}=0.
\end{equation}

The following result is, obviously, important
\begin{prop} \label{propos12}  If $N$ is either an a.P.c.
or a quasi-Dirac submanifold, the bivector field
$\Pi'=p_T(\Pi|_N)$ is a Poisson bivector field on $N$. Moreover,
in the a.P.c. case $\Pi'$ does not depend on the choice of the
normal bundle among those which satisfy Definition \ref{algebr}.
\end{prop} \noindent{\bf Proof.} From (\ref{Pilocal}), it follows
\begin{equation} \label{Pi'}
\Pi'=\frac{1}{2}\left(\Pi^{uv}\frac{\partial}{\partial y^u}
\wedge\frac{\partial}{\partial y^v}\right)_{x^c=0}.\end{equation}
Then, from $[\Pi,\Pi]=0$ and either (\ref{Piauzero}) or
(\ref{derivPi}) we get
\begin{equation} \label{uvw} [\Pi,\Pi]^{u_1u_2u_3}|_{x^c=0}=
2\left(
\sum_{Cycl(u_1,u_2,u_3)}\Pi^{u_1w}\frac{\partial\Pi^{u_2u_3}}{\partial
y^w}\right)_{x^c=0}\end{equation} $$=[\Pi',\Pi']^{u_1u_2u_3}=0$$
i.e., $\Pi'$ is a Poisson bivector field on $N$.

Finally, (\ref{Piauzero}) shows that
$\iota_*(\sharp_{\Pi'}(dy^u))=\sharp_{\Pi}(dy^u)$
$(\iota:N\subseteq M)$, and this proves the last assertion.
Q.e.d.

The Poisson structure $\Pi'$ is said to be {\it induced} by the
Poisson structure $\Pi$, and was defined and studied for Dirac
submanifolds in \cite{Xu}. In the a.P.c. case, the submanifold has
a {\it second fundamental form} $(L_X\Pi)|_{ann(\nu N)}$, which
vanishes iff $N$ is a Dirac submanifold.

If $\sigma:W\rightarrow N$ is a compatible tubular neighborhood of
$(N,\nu N)$, the induced Poisson structure is characterized by
\begin{equation}\label{caracttub}
\{f,g\}_{\Pi'}= \{f\circ\sigma,g\circ\sigma\}_{\Pi}\circ\iota,
\end{equation} $\forall f,g\in C^\infty(N)$ and $\iota:N\subseteq M$.
With the same tubular neighborhood, the a.P.c. property is
equivalent with the fact that $\forall f\in C^\infty(N)$, the
$\Pi$-hamiltonian vector field of $f\circ\sigma$ is tangent to $N$
or that one has \begin{equation}\label{compD}
\iota_*\circ\sharp_{\Pi'}=\sharp_{\Pi}\circ\sigma^{*}.
\end{equation} Furthermore,
$N$ is a Dirac submanifold if, besides the above, it also has the
property that $\forall x\in N$, $\forall X\in T_x(W_x)$, $\forall
f,g\in C^\infty(N)$ one has
\begin{equation}\label{ptDirac} X\{f\circ\sigma,g\circ\sigma\}_\Pi=0.
\end{equation}

We summarize the above remarks in \begin{prop}\label{propos13} A
submanifold $N$ of the Poisson manifold $(M,\Pi)$ is an a.P.c.
submanifold iff $N$ is endowed with a Poisson structure $\Pi'$ and
has a tubular neighborhood $\sigma:W\rightarrow N$ such that
conditions (\ref{caracttub}) and (\ref{compD}) hold. Furthermore,
$N$ is a Dirac submanifold iff it is a.P.c. and condition
(\ref{ptDirac}) holds too. \end{prop}
\begin{rem}\label{obs0} {\rm In \cite{Xu} the author uses the independence
of the induced Poisson structure on the choice of the normal
bundle to define {\it local Dirac submanifolds} \cite{Xu} as
submanifolds $N$ of the Poisson manifold $(M,\Pi)$ such that,
$\forall x\in N$, there exists an open neighborhood $U$ in $M$
where $N\cap U$ is a Dirac submanifold. Proposition \ref{propos12}
shows that a local Dirac submanifold also inherits a well defined,
global, induced Poisson structure. In fact, Proposition
\ref{propos12} shows that {\it local a.P.c. submanifolds} may be
defined similarly.}
\end{rem}
\begin{rem}\label{obs1} {\rm Proposition
\ref{propos13} suggests considering  submanifolds $N$ of $(M,\Pi)$
which come endowed with a Poisson structure $\Pi'$ such that, for
some tubular neighborhood $\sigma:W\rightarrow N$, $\sigma$ is a
Poisson mapping or, equivalently, the brackets
$\{f\circ\sigma,g\circ\sigma\}|_{\Pi}$, which are defined on $W$,
are constant along the fibers of $\sigma$. Such submanifolds
deserve the name of {\it strong Dirac submanifolds}. Obviously,
they are Dirac submanifolds, and, with respect to adapted
coordinates, one must have $\partial\Pi^{uv}/\partial
x^a\equiv0$.}
\end{rem}

The complete lift of the Poisson bivector field $\Pi$ is a Poisson
structure $\Pi^C$ called the {\it tangent Poisson structure} of
$\Pi$. The tangent structure is exactly the one induced  by the
Lie algebroid structure of $T^*M$ defined by $\Pi$ on the total
space of its dual vector bundle $TM$. From Proposition
\ref{propos01}, we get the following characteristic property of
Dirac submanifolds
\begin{prop}\label{propos15} {\rm\cite{Xu}} The submanifold $N$ of $(M,\Pi)$ is
a Dirac submanifold iff there exists a normal bundle $\nu N$ which
is a coisotropic submanifold of $(TM,\Pi^C)$.\end{prop}

The following result is significant for the next section. Recall
that a Poisson structure $\Pi$ is {\it homogeneous} if there
exists a vector field $Z\in\Gamma TM$ such that
\begin{equation} \label{homog} L_Z\Pi=-\Pi. \end{equation}
\begin{prop}\label{propos19} Let $N$ be an a.P.c. submanifold of the
homogeneous Poisson manifold $(M,\Pi,Z)$ such that $Z|_N\in \Gamma
TN$. Then $(N,\Pi',Z|_N)$, where $\Pi'$ is the induced Poison
structure, also is a homogeneous Poisson manifold.
\end{prop} \noindent{\bf Proof.} If the homogeneity
condition (\ref{homog}) is evaluated on $(d\varphi,d\psi)$
$(\varphi,\psi\in C^\infty(M))$, one gets the equivalent condition
\begin{equation} \label{homog2} \{\varphi,\psi\}_{\Pi}=
Z\{\varphi,\psi\}_{\Pi}-[Z,X^{\Pi}_\varphi]\psi+[Z,X^{\Pi}_\psi]\varphi,
\end{equation} where $X$ denotes hamiltonian vector fields.

Now, let $\sigma:W\rightarrow N$ be a tubular neighborhood of $N$
where the conditions of Proposition \ref{propos13} hold. Then,
(\ref{homog2}), written for
$\varphi=f\circ\sigma,\psi=g\circ\sigma$ $(f,g\in C^\infty(M))$
and composed by $\iota$, provides the similar condition
(\ref{homog2}) for $\Pi'$. Q.e.d.

Now, we point out the existence of a specific cohomology related
with a Dirac submanifold $(N,\nu N)$ of the Poisson manifold
$(M,\Pi)$. Namely, since the Poisson bivector field $\Pi$ is
soldered to $N$,
\begin{equation} \label{newcoh} \mathcal{V}(M,N,\nu N)=
(\mathcal{V}^k(M,N,\nu N),\;\partial_{\Pi}=-[\Pi,\,.])
\end{equation} is a subcomplex of the Lichnerowicz-Poisson cochain
complex of $(M,\Pi)$ (e.g., \cite{V1}), therefore, it defines
cohomology spaces $H^k_{sP}(M,N,\nu N)$ that will be called {\it
soldered Poisson cohomology spaces}.
\begin{prop} \label{newprop} The homomorphism $\sharp_{\Pi}:
T^*M\rightarrow TM$ induces homomorphisms
\begin{equation} \label{homocoh} \sharp_{\Pi}^k:
H^k_{sde R}(M,N,\nu N)\longrightarrow H^k_{sP}(M,N,\nu N).
\end{equation} If the Poisson structure $\Pi$ is defined by a
symplectic form $\omega$, the mappings (\ref{homocoh}) are
isomorphisms.
\end{prop}
\noindent{\bf Proof.} Using adapted local coordinates, it is easy
to check that, $\forall\kappa\in \Omega^k(M,N,\nu N)$,
$\sharp_{\Pi}\kappa\in\mathcal{V}^k(M,N,\nu N)$ and, as for the
general Poisson cohomology,
$\partial_{\Pi}\circ\sharp_{\Pi}=-\sharp_{\Pi}\circ d$. This
justifies the existence of  the homomorphisms (\ref{homocoh}).

In the symplectic case, $N$ must be a symplectic
$(2n)$-dimensional submanifold of $(M^{2m},\omega)$, and $\nu N$
must be the $\omega$-orthogonal bundle of $TN$ \cite{Xu}. Indeed,
each point $x\in N$ has an open neighborhood with coordinates
$(x^a,x^{a^*},y^u,y^{u^*})$
$(a=1,...,m-n,\,u=2m-2n+1,...,2m-n,\,a^*=a+m-n,\,u^*=u+n)$ such
that
\begin{equation} \label{omegacan}
\omega=\sum_adx^a\wedge dx^{a^*}+ \sum_udy^u\wedge
dy^{u^*},\end{equation} and $N$ has the local equations
$x^a=0,x^{a^*}=0$ \cite{Mar}. Obviously, these coordinates also
are adapted coordinates with respect to the $\omega$-orthogonal
bundle $\nu N$ of $TN$. The uniqueness of the normal bundle
follows from (\ref{cond1theta}). Then, $\sharp_{\Pi}$ has the
inverse $-\flat_\omega$, and the latter induces the inverses of
the homomorphisms (\ref{homocoh}). Q.e.d.
\section{Dirac submanifolds of Jacobi manifolds}
For a detailed study of Jacobi manifolds we refer the reader to
\cite{DML} and its references. Jacobi manifolds are a natural
generalization of Poisson manifolds namely, a Jacobi structure on
a manifold $M^m$ is a Lie algebra bracket $\{f,g\}$ on
$C^\infty(M)$, which is given by bidifferential operators. It
follows that one must have
\begin{equation} \label{Jacobi} \{f,g\}=\Lambda(df,dg)+fEg-gEf,
\end{equation} where $E$ is a vector field and $\Lambda$ is a
bivector field on $M$ such that\footnote{The minus sign comes from
our sign convention for the Schouten-Nijenhuis bracket in Section
1.}
\begin{equation} \label{Jacobi2}
[\Lambda,\Lambda]=-2E\wedge\Lambda,\;L_E\Lambda=0.\end{equation}
Thus, if $E=0$ we have a Poisson structure.

The Jacobi structure $(M,\Lambda,E)$ is equivalent with the
homogeneous Poisson structure \begin{equation} \label{Pois-Jacobi}
\Pi=e^{-\tau}(\Lambda+\frac{\partial}{\partial\tau}\wedge E),\;\;
Z=\frac{\partial}{\partial\tau} \end{equation} on
$M\times\mathbb{R}$; $M$ will then be identified with
$M\times\{0\}$. For instance \cite{Lich}, let $\Pi$ be a linear
Poisson structure on $ \mathbb{R}^n\backslash\{0\}$, and consider
the diffeomorphism $
S^{n-1}\times\mathbb{R}\approx\mathbb{R}^n\backslash\{0\}$ defined
by $$x^i=e^{-\tau}u^i, \;\;\sum_{i=1}^n (u^i)^2=1,\;\;
\tau\in\mathbb{R}.$$ Then, it is easy to check that $\Pi$ must be
of the form (\ref{Pois-Jacobi}), which provides a Jacobi structure
on $S^{n-1}$. We call this structure a {\it Lichnerowicz-Jacobi
structure} of the sphere.

On a Jacobi manifold, one may define {\it hamiltonian vector
fields}
\begin{equation} \label{ham-Jacobi} X_f=\sharp_\Lambda
df+fE\;\;(f\in C^\infty(M)),\end{equation} and they span a
generalized foliation $ \mathcal{S}$ such that the leaves of $
\mathcal{S}$ are either contact or locally conformal symplectic
manifolds. For instance, in the case of a Lichnerowicz-Jacobi
structure the leaves are the orbits of the {\it quotient coadjoint
representation} of a connected Lie group $G$ with the Lie algebra
$ \mathcal{G}$ of structure defined by the corresponding linear
Poisson structure $\Pi$ (see above), i.e., the action defined on
$S^{n-1}$ by the coadjoint action of $G$ on $
\mathcal{G}^*\backslash\{0\}\approx \mathbb{R}^n\backslash\{0\}$
if $S^{n-1}$ is seen as a quotient space of $
\mathcal{G}^*\backslash\{0\}$) \cite{Lich}.

Another important fact we want to recall is that, $\forall
\varphi\in C^\infty(M)$, the bracket
\begin{equation}\label{conform1} \{f,g\}^\varphi=e^{-\varphi}
\{e^\varphi f,e^\varphi g\} \end{equation} is a Jacobi bracket
said to have been obtained by a {\it conformal change} of the
original bracket. The tensor fields of the new bracket are
\begin{equation} \label{conform2}
\Lambda^\varphi=e^\varphi\Lambda,\;E^\varphi=e^\varphi(E+i(d\varphi)\Lambda).
\end{equation}

Now, we begin our considerations on submanifolds.
\begin{defin} \label{Dir-Jacobi} {\rm Let $(N,\nu N)$ be a normalized
submanifold of the Jacobi manifold $(M,\Lambda,E)$. Then: 1) $N$
is an {\it almost Dirac submanifold} if $\Lambda$ is an
$N$-soldered bivector field; 2) $N$ is an {\it algebraically
Jacobi-compatible (a.J.c.)} submanifold if $\Lambda$ and $E$ are
algebraically compatible with the normalization of $N$; 3) $N$ is
a {\it (quasi-) Dirac submanifold} if $\Lambda$ and $E$ are
(quasi-) soldered to $(N,\nu N)$.}
\end{defin}

Equivalently, $N$ is an almost Dirac submanifold if each point
$x\in N$ has a neighborhood with adapted coordinates as in Section
1, such that
\begin{equation} \label{Lambdalocal}
\Lambda=\frac{1}{2}\Lambda^{ab}\frac{\partial}{\partial
x^a}\wedge\frac{\partial}{\partial x^b} +
\Lambda^{au}\frac{\partial}{\partial
x^a}\wedge\frac{\partial}{\partial y^u} +
\frac{1}{2}\Lambda^{uv}\frac{\partial}{\partial
y^u}\wedge\frac{\partial}{\partial y^v}, \end{equation} where
\begin{equation} \label{Lauzero} \Lambda^{au}|_{x^c=0}=0, \end{equation}
\begin{equation} \label{derivL} \left.\frac{\partial\Lambda^{uv}}{\partial
x^a}\right|_{x^c=0}=0. \end{equation} Then, $N$ is a Dirac
submanifold if, furthermore, the vector field $E$ is tangent to
$N$ and has the local expression \begin{equation} \label{Elocal}
E=\epsilon^a\frac{\partial}{\partial
x^a}+\epsilon^u\frac{\partial}{\partial y^u},\end{equation} where
\begin{equation}\label{condE} \epsilon^a|_{x^c=0}=0,\hspace{5mm}
\left.\frac{\partial\epsilon^u}{\partial x^a}\right|_{x^c=0}=0.
\end{equation}
For the quasi-Dirac case, we only have the condition
(\ref{derivL}) and the second equality (\ref{condE}). Finally, $N$
is an a.J.c. submanifold if (\ref{Lauzero}) and the first
condition (\ref{condE}) hold.
\begin{prop} \label{propos21} Let $(N,\nu N)$ be either an almost
Dirac or an a.J.c. or a quasi-Dirac submanifold of the Jacobi
manifold $(M,\Lambda,E)$. Then,
$[\Lambda'=p_T(\Lambda|_N),E'=p_T(E|_N)]$ is a Jacobi structure on
$N$. Furthermore, in the almost Dirac and the a.J.c. case
$\Lambda'$ does not depend on the choice of the normalization.
\end{prop} \noindent{\bf Proof.} With local adapted coordinates we
have \begin{equation} \label{Lambda'}
\Lambda'=\frac{1}{2}\left(\Lambda^{uv}\frac{\partial}{\partial
y^u} \wedge\frac{\partial}{\partial
y^v}\right)_{x^c=0},\end{equation} and
\begin{equation} \label{uvwLambda}
[\Lambda,\Lambda]^{u_1u_2u_3}|_{x^c=0}=
2\left(\sum_{Cycl(u_1,u_2,u_3)}\Lambda^{u_1a}
\frac{\partial\Lambda^{u_2u_3}}{\partial x^a}\right.\end{equation}
$$+\left.\sum_{Cycl(u_1,u_2,u_3)}
\Lambda^{u_1w}\frac{\partial\Lambda^{u_2u_3}}{\partial
y^w}\right)_{x^c=0}=[\Lambda',\Lambda']^{u_1u_2u_3}.$$ Hence, in
all the cases of the proposition we have
$$[\Lambda',\Lambda']=p_T([\Lambda,\Lambda]|_N)=-2E'\wedge\Lambda'.$$

Then, an examination of the coordinate expression of $L_E\Lambda$,
where $E$ and $\Lambda$ are given by (\ref{Elocal}) and
(\ref{Lambdalocal}), respectively, shows that the conditions
({\ref{Lauzero}) and (\ref{derivL}), as well as either
(\ref{Lauzero}) and the first condition (\ref{condE}) or
(\ref{derivL}) and the second condition (\ref{condE}), imply
$$L_{E'}\Lambda'=p_T(L_E\Lambda|_N)=0.$$

Finally, where asserted, the independence of $\Lambda'$ of the
normalization follows from
$\iota_*(\sharp_{\Lambda'}(dy^u)=\sharp_\Lambda(dy^u)$,
$\iota:N\subseteq M$. Q.e.d.

We notice that formula (\ref{uvwLambda}) also implies
\begin{prop} \label{corolar21}  Let $(N,\nu N)$ be a normalized
submanifold of the Jacobi manifold $(M,\Lambda,E)$ such that
$\Lambda$ is algebraically compatible with the normalization and
$E$ is a normal field of $N$. Then $\Lambda'$ is a Poisson
structure on $N$, and it is independent on the choice of $\nu N$
among all possible choices that contain $E|_N$.
\end{prop}

The Jacobi or Poisson structures defined on $N$ by $(\Lambda',E')$
are said to be {\it induced} by $(\Lambda,E)$. In the cases where
only algebraic compatibility holds, invariants of the second
fundamental form type $(L_X\Lambda)|_{ann(\nu N)}$,
$[X,E]|_{ann(\nu N)}$, where $X$ is normal to $N$, appear.

Proposition \ref{proposinvol} allows us to give some simple
examples. Consider the Jacobi manifold $M=\mathbb{R}^{3n+1}$ with
\begin{equation} \label{exampleeq} \Lambda= \sum_{i=1}^nu_i
\frac{\partial}{\partial q^i}\wedge\frac{\partial}{\partial p_i}
+(t\frac{\partial}{\partial
t})\wedge(\sum_{j=1}^np_j\frac{\partial}{\partial p_j}),\;
E=t\frac{\partial}{\partial t} \end{equation} (the variables of
(\ref{exampleeq}) are the natural coordinates of $M$). Then the
hyperplane $t=0$ is the fixed point locus of the involution
$(u_i,q^i,p_i,t)\rightarrow(u_i,q^i,p_i,-t)$. This involution
preserves the tensor fields (\ref{exampleeq}) hence, the
hyperplane $t=0$ is a Dirac submanifold with an induced Poisson
structure. For the same $M$, the involution
$(u_i,q^i,p_i,t)\rightarrow(-u_i,-q^i,p_i,t)$ also preserves
(\ref{exampleeq}) hence, its fixed point locus, which is the
$(n+1)$-plane $u_i=0,q^i=0$, is a Dirac submanifold with an
induced Jacobi structure. Finally, if we restrict $M$ to the
domain $p_i>0,t>0$ and consider the involution that sends $t$ to
$1/t$ and preserves the other coordinates, only $\Lambda$ of
(\ref{exampleeq}) is preserved hence, the fixed point locus, which
is the hyperplane $t=1$, is an almost Dirac submanifold. Moreover,
the last involution sends $E$ to $-E$, therefore, $E$ is normal to
the submanifold, and the induced structure is a Poisson structure.

Another interesting fact is \begin{prop} \label{propos22} Let
$(M,\Lambda,E)$ be a Jacobi manifold, and $(M\times
\mathbb{R},\Pi)$, with $\Pi$ defined by {\rm(\ref{Pois-Jacobi})},
the corresponding homogeneous Poisson manifold. Then, $N$ is an
a.J.c. or a Dirac submanifold of the former iff $N\times
\mathbb{R}$ is an a.J.c., respectively Dirac, submanifold of the
latter.
\end{prop} \noindent{\bf Proof.} We will use the lift of $\nu N$
to $N\times \mathbb{R}$ as a normal bundle. $\tau$ of
(\ref{Pois-Jacobi}) is a coordinate along $N\times \mathbb{R}$,
and we see that if $\Lambda$ and $E$ are of the local form
(\ref{Lambdalocal}), (\ref{Elocal}) then $\Pi$ satisfies the
conditions (\ref{Piauzero}), (\ref{derivPi}), and conversely.
Q.e.d.

From Proposition \ref{propos01} it follows that the almost Dirac
and Dirac submanifolds of a Jacobi manifold may also be
characterized by using the tangent bundle $TM$ namely,
\begin{prop} \label{propos23} The normalized submanifold $(N,\nu
N)$ is almost Dirac iff the submanifold $\nu N$ of $TM$ is such
that $\sharp_{\Lambda ^C}(ann(T\nu N))\subseteq T\nu N$, where
$\Lambda^C$ is the complete lift of $\Lambda$. Furthermore, $N$ is
a Dirac submanifold iff besides the previous condition, one also
has $E^C\in\Gamma(T\nu N)$, where $E^C$ is the complete lift of
the vector field $E$. \end{prop}
\begin{rem} \label{obs21} {\rm
The tensor fields $(\Lambda^C,E^C)$ do not define a Jacobi
structure on $TM$. A {\it tangent Jacobi structure} can be
obtained by considering the Poisson structure induced on the
manifold $TM\times \mathbb{R}$ by the Lie algebroid
$J^1M=T^*M\times \mathbb{R}$  of $(\Lambda,E)$ \cite{Kerb}.
Namely, with the notation of (\ref{completeform}),
(\ref{completemult}), if we associate with each cross section
$(f,\alpha_idz^i)\in\Gamma J^1M$ the function
$e^\tau(f+\alpha_i\hspace{-1mm}\stackrel{.}{z}^i)$ $\in
C^\infty(TM\times\mathbb{R})$, the Lie algebroid bracket of $J^1M$
yields a Poisson bracket of the specified kind of functions, which
extends to a Poisson bracket on $C^\infty(TM\times\mathbb{R})$.
Computations show that the Poisson bivector of this structure is
\begin{equation} \label{tangentDir}
\tilde\Pi=e^{-\tau}[\Lambda^C-\Lambda^V-\mathcal{E}
\wedge(E^C-E^V)+\frac{\partial}{\partial \tau}\wedge E^C],
\end{equation} where the upper index $V$ denotes the {\it vertical
lift} \cite{YI}, and $ \mathcal{E}$ is the {\it Euler vector
field} i.e., \begin{equation} \label{vertical} E^V= \stackrel{.}{
z}^i\frac{\partial} {\partial\hspace{-1mm}\stackrel{.}{ z}^i},
\;\;\Lambda^V=\frac{1}{2}\Lambda^{ij}\frac{\partial}
{\partial\hspace{-1mm}\stackrel{.}{ z}^i}\wedge\frac{\partial}
{\partial\hspace{-1mm}\stackrel{.}{ z}^j},\;\; \mathcal{E}=
\stackrel{.}{ z}^i\frac{\partial}
{\partial\hspace{-1mm}\stackrel{.}{ z}^i}.\end{equation}
Accordingly $(\Lambda^C-\Lambda^V-\mathcal{E} \wedge(E^C-E^V),
E^C)$ is a Jacobi structure on $TM$, which deserves the name of
{\it tangent Jacobi structure}.}\end{rem}

A particular class of Dirac submanifolds was studied in
\cite{DML}, and we reprove here
\begin{prop} \label{propos24} {\rm \cite{DML}} Assume that $N$ is a submanifold of
$(M,\Lambda,E)$ such that $\sharp_\Lambda(ann(TN))$ is a normal
bundle $\nu N$ of $N$. Then $N$ is a Dirac submanifold iff the
vector field $E$ is tangent to $N$. Furthermore, there always
exists a conformal change (\ref{conform1}), with $\varphi|_N=0$,
such that $N$ is a Dirac submanifold of
$(M,\Lambda^\varphi,E^\varphi)$.
\end{prop} \noindent{\bf Proof.} If $N$ is a Dirac submanifold,
$E\in\Gamma TN$ by definition. For the converse, we use the normal
bundle of the hypothesis, and represent $\Lambda$ and $E$ by
(\ref{Lambdalocal}) and (\ref{Elocal}), respectively. Clearly, the
choice of $\nu N$ is such that $\sharp_\Lambda(dx^a)\in\Gamma\nu
N$ along $N$, which is equivalent to (\ref{Lauzero}).  If this
condition holds, the ${auv}$-component of the first equality
(\ref{Jacobi2})yields
\begin{equation} \label{Lambdaauv}
(\Lambda^{ab}\frac{\partial\Lambda^{uv}}{\partial
x^b}+E^a\Lambda^{uv})|_N=0. \end{equation} Since
$\sharp_\Lambda(ann(TN))$ is normal to $N$ iff the matrix
$(\Lambda^{ab})$ is non degenerate, we see that $E$ tangent to $N$
implies (\ref{derivL}). Furthermore, if $E|_N\in\Gamma TN$, the
second equality (\ref{Jacobi2}) yields
\begin{equation}\label{Eu} \left.\Lambda^{ab}\frac{\partial
E^u}{\partial x^b}\right|_N=0,
\end{equation} therefore (\ref{condE}) holds.

Concerning the last part of the proposition, (\ref{conform2})
shows that the required conformal transformation exists if there
exists a function $\varphi\in C^\infty(M)$, which vanishes on $N$
and is such that \begin{equation} \label{eqptphi}
E^a|_N=\left.\left(\Lambda^{ab}\frac{\partial\varphi}{\partial
x^b}\right)\right|_N. \end{equation} Since $(\Lambda^{ab})$ is non
degenerate, the conditions for $\varphi$ prescribe the $1$-jet
with respect to the variables $(x^a)$ of $\varphi$ at the points
of $N$. Therefore, a required function $\varphi$ exists around
every point of $N$. Then, these local solutions may be glued up by
a partition of unity along $N$. (See also the argument of
\cite{DML}). Q.e.d.

We end by a discussion of Dirac submanifolds of the {\it
transitive} Jacobi manifolds i.e., locally conformal symplectic
(l.c.s.) and contact manifolds \cite{DML}.
\begin{prop} \label{propos31} A submanifold $N$ is an almost Dirac
submanifold of the l.c.s. manifold $M$ iff it is Dirac, and this
happens iff $N$ inherits from $M$ an induced l.c.s. structure.
Moreover, there is only one possible normal bundle, the symplectic
orthogonal bundle of $TN$. \end{prop} \noindent{\bf Proof.} Recall
that the l.c.s. structure of $M$ is a non degenerate $2$-form
$\Omega$ such that for some open covering $M=\cup_\alpha
U_\alpha$, $\forall\alpha$,
$\Omega|_{U_\alpha}=e^{-\sigma_\alpha}\Omega_\alpha$, where
$\sigma_\alpha$ are functions, $\Omega_\alpha$ are $2$-forms and
$d\Omega_\alpha=0$. Equivalently, $d\Omega=\omega\wedge\Omega$,
where $\omega$ is the closed $1$-form defined by gluing up the
local forms $d\sigma_\alpha$ ($\omega$ is called the {\it Lee
form}). It is known that $M$ is a Jacobi manifold with the
structure defined by the bivector field $\Lambda$, where
$\sharp_\Lambda=\flat_\Omega^{-1}$, and the vector field
$E=\sharp_\Lambda\omega$ \cite{DML}.

Assume that $N$ is an almost Dirac submanifold with the normal
bundle $\nu N$. Then, $\sharp_\Lambda|_{ann(TN)}$ is an
isomorphism onto $\nu N$, which, just like (\ref{cond1theta}),
ensures the uniqueness of $\nu N$, and $\sharp_\Lambda|_{ann(\nu
N)}$ is an isomorphism on $TN$, which is equivalent with the fact
that $\iota^*\Omega$ $(\iota:N\subseteq M)$ is non degenerate and
provides an l.c.s. structure on $N$. Accordingly, we may use again
Marle's theorem \cite{Mar}, and find local coordinates of $N$ on
some neighborhood $U_\alpha$ such that, with the notation of
(\ref{omegacan}), one has
\begin{equation} \label{canoniclcs} \Omega|_{U_\alpha}=
e^{-\sigma_\alpha}(\sum_adx^a\wedge dx^{a^*} + \sum_udy^u\wedge
dy^{u^*}). \end{equation} Obviously, this expression of $\Omega$
implies that $\nu N$ is $\Omega$-orthogonal to $TN$ and that the
coordinates used in (\ref{canoniclcs}) are adapted coordinates.

Now, condition(\ref{derivL}) applied to (\ref{canoniclcs}) becomes
$(\partial\sigma_\alpha/\partial x^a)_{x^c=0}=0$ i.e.,
$\omega_a|_{x^c=0}$ $=0$. Furthermore, one of the conditions that
express $d\omega=0$ is $\partial\omega_a/\partial y^u=
\partial\omega_u/\partial x^a$, whence we also get
$(\partial\omega_u/\partial x^a)_{x^c=0}=0$. Therefore, the Lee
form $\omega$, and the vector field $E$ too, are soldered to $N$,
and $N$ must be a Dirac submanifold of $M$. The converse part of
the proposition follows from (\ref{canoniclcs}). Q.e.d.
\begin{prop} \label{contactD} Let $M^{2m+1}$ be a contact manifold
with the contact $1$-form $\theta$. Then a submanifold
$\iota:N\subseteq M$ is a Dirac submanifold iff $\iota^*\theta$ is
a contact form on $N$. Furthermore, the normal bundle of $N$ is
unique, and it is the $d\theta$-orthogonal bundle of $TN$.
\end{prop} \noindent{\bf Proof.} Recall that $\theta$ is a contact
form iff $\theta\wedge (d\theta)^m$ vanishes nowhere. A contact
form produces a Jacobi structure \cite{DML}, which consists of the
vector field $E$  defined by
\begin{equation} \label{Reeb} i(E)\theta=1,\;i(E)d\theta=0,
\end{equation} and the bivector field
\begin{equation} \label{Lambdacontact}
\Lambda(df,dg)=d\theta(X^\theta_f,X^\theta_g)\;\;(f,g\in
C^\infty(M)), \end{equation} where the {\it hamiltonian vector
field} $X^\theta_f$ is defined by \begin{equation}
\label{hamiltcontact}
i(X^\theta_f)\theta=f,\;i(X^\theta_f)d\theta=-df+(Ef)\theta.
\end{equation} From (\ref{hamiltcontact}) we get
\begin{equation} \label{Lambdanou}X^\theta_f=\sharp_\Lambda
df+fE,\;\Lambda(df,dg)=d\theta(\sharp_\Lambda df,\sharp_\Lambda
dg).\end{equation}

If $M$ is the contact manifold above, $M\times \mathbb{R}$ has the
Poisson bivector $\Pi$ given by (\ref{Pois-Jacobi}), and it also
has the symplectic form \begin{equation} \label{symplprodus}
\Omega=e^\tau(d\theta+d\tau \wedge\theta). \end{equation} An easy
computation shows that all the functions of the form $e^\tau f\in
C^\infty(M\times \mathbb{R})$ $(\tau\in\mathbb{R},f\in
C^\infty(M))$ have the same hamiltonian vector fields with respect
to $\Pi$ and $\Omega$, therefore,
$\sharp_\Pi\circ\flat_\Omega=-Id.$

Now, Proposition \ref{propos22} tells that $N$ is a Dirac
submanifold of $M$ iff $N\times \mathbb{R}$ is a Dirac submanifold
of $M\times \mathbb{R}$. In the present case, this means that
$N\times \mathbb{R}$ is a symplectic submanifold of $(M\times
\mathbb{R},\Omega)$, and it follows that $(N,\iota^*\theta)$ must
be a contact manifold, and that the normal bundle must be the one
indicated by the proposition.  Q.e.d.
\small{} 
\hspace*{7.5cm}{\small \begin{tabular}{l} Department of
Mathematics\\ University of Haifa, Israel\\ E-mail:
vaisman@math.haifa.ac.il \end{tabular}}
\end{document}